\documentclass[11pt,twoside]{article}
\usepackage{geometry}                		
\usepackage{graphicx}									
\usepackage{amssymb}\vspace{0.25cm}
\usepackage{mathtools}
\usepackage{setspace}
\usepackage{tikz-cd}
\usepackage[mathscr]{euscript}
\newcommand{\abs}[1]{\left\vert#1\right\vert}
\usepackage[normalem]{ulem}
\usepackage{manfnt}
\usepackage{pgfplots}
\usepackage{pifont}
\usepackage[safe]{tipa}
\usepackage{marvosym}
\usepackage{scalerel,stackengine}
\usepackage{titlesec}
\usepackage{makeidx}
\usepackage{centernot}

\usepackage{xspace}
\usepackage[noauto]{chappg}
\usepackage[OT2,T1]{fontenc}
\usepackage{hyperref}
\usepackage{amsmath}
\usepackage[amsmath,thmmarks]{ntheorem}

\usepackage{tabularx}
\newcolumntype{Y}{>{\raggedright\arraybackslash}X}

\newcommand\mA{%
$A$\xspace
}

\newcommand\mG{%
$G$\xspace
}

\newcommand\mK{%
$K$\xspace
}

\newcommand\R{%
\mathbb{R}
}

\newcommand\Z{%
\mathbb{Z}
}

\newcommand\bc{%
\textbf{c}\xspace
}

\newcommand\bCon{%
\textbf{Con}\xspace
}
\newcommand\bSp{%
\textbf{Sp}\xspace
}

\newcommand\bA{%
\textbf{A}\xspace
}
\newcommand\bB{%
\textbf{B}\xspace
}
\newcommand\bG{%
\textbf{G}\xspace
}

\newcommand\bInd{%
\textbf{Ind}\xspace
}

\newcommand\bGL{%
\textbf{GL}\xspace
}
\newcommand\bSL{%
\textbf{SL}\xspace
}

\newcommand\sC{%
\mathcal{C}
}

\newcommand\sO{%
\mathcal{O}
}

\newcommand\rank{%
\text{rank}\hspace{0.05cm}
}

\newcommand\con{%
\text{con}\hspace{0.05cm}
}

\newcommand\dis{%
\text{dis}\hspace{0.05cm}
}
\newcommand\Ker{%
\text{Ker}\hspace{0.05cm}
}
\newcommand\tr{%
\text{tr}\hspace{0.05cm}
}

\newcommand\restr[2]{%
{#1}|{#2}
}

\newcommand{\norm}[1]{\left\lVert #1 \right\rVert}

\usepackage{scalerel,stackengine}
\stackMath
\newcommand\reallywidehat[1]{%
\savestack{\tmpbox}{\stretchto{%
  \scaleto{%
    \scalerel*[\widthof{\ensuremath{#1}}]{\kern-.6pt\bigwedge\kern-.6pt}%
    {\rule[-\textheight/2]{1ex}{\textheight}}
  }{\textheight}%
}{0.5ex}}%
\stackon[1pt]{#1}{\tmpbox}%
}
\DeclareFontFamily{U}{wncy}{}
\DeclareFontShape{U}{wncy}{m}{n}{<->wncyr10}{}
\DeclareSymbolFont{mcy}{U}{wncy}{m}{n}
\DeclareMathSymbol{\Sh}{\mathord}{mcy}{"58} 

\makeatletter
\DeclareRobustCommand\widecheck[1]{{\mathpalette\@widecheck{#1}}}
\def\@widecheck#1#2{%
    \setbox\z@\hbox{\m@th$#1#2$}%
    \setbox\tw@\hbox{\m@th$#1%
       \widehat{%
          \vrule\@width\z@\@height\ht\z@
          \vrule\@height\z@\@width\wd\z@}$}%
    \dp\tw@-\ht\z@
    \@tempdima\ht\z@ \advance\@tempdima2\ht\tw@ \divide\@tempdima\thr@@
    \setbox\tw@\hbox{%
       \raise\@tempdima\hbox{\scalebox{1}[-1]{\lower\@tempdima\box
\tw@}}}%
    {\ooalign{\box\tw@ \cr \box\z@}}}
\makeatother

\newtheoremstyle{xx}
  {4pt}
  {0pt}
  {\upshape}
  {\bfseries}
  {}
  { }
  {}
  
\makeatletter 
 \newtheoremstyle{myu}%
  {\upshape\item[ \indent\indent\bf\underline{\theorem@headerfont ##2:}]}%
\makeatother
\makeatletter 
 \newtheoremstyle{myn}%
  {\item[\hskip\labelsep \ \bf ##1 \theorem@headerfont ##2.]}%
\makeatother
\theoremstyle{myn}
\theoremstyle{myu}
{\upshape}

\titleformat{\chapter}[display]
{\normalfont\filcenter\huge\bfseries}{}{0pt}{\large}

\usepackage[OT2,OT1]{fontenc} 
\newcommand\cyr
{
\renewcommand\rmdefault{wncyr} 
\renewcommand\sfdefault{wncyss} 
\renewcommand\encodingdefault{OT2} 
\normalfont
\selectfont
}

\DeclareTextFontCommand{\textcyr}{\cyr}

\makeindex 
\begin{document}

\title{\textbf{EXPLICIT FORMULAS FROM THE\\
CONTINUOUS SPECTRUM}}
\vspace{1cm}
\author{
\vspace{.3cm}\\
by\vspace{.3cm}\\ 
\vspace{.3cm}
M. Scott Osborne$^{\ *}$\\
\vspace{.3cm}
and\\
\vspace{.3cm}
Garth Warner
\footnote{ \ Research supported in part by the National Science Foundation.}}


\date{\vspace{.3cm} University of Washington\\
Seattle, Washington 98195}

\maketitle                              

\titlespacing*{\chapter}{0pt}{-50pt}{40pt}
\setlength{\parskip}{0.1em}
\include{_tocX}
\include{_Preface}


\[
\textbf{EXPLICIT FORMULAS FROM THE CONTINUOUS SPECTRUM}
\]
\[
\text{Scott Osborne and Garth Warner}
\]

\setlength\parindent{2em}
\setcounter{theoremn}{0}
\vspace{.5cm}

{\small
The purpose of this note is to announce the results of our investigation into the role played by the continuous spectrum in the development of the Selberg trace formula vis-\`a-vis a pair $(G,\Gamma)$.  
For the sake of simplicity, we shall restrict ourselves to a ``rank-2'' situation, a case in point being when
\[
\begin{cases}
\ G = \bSL(3,\R)\\
\ \Gamma = \bSL(3,\Z)
\end{cases}
.
\]
Full details (in all generality) will appear elsewhere.\\
}
\vspace{1.cm}

Let \mG be a reductive Lie group, $\Gamma$ a lattice in \mG, both subject to the usual conditions (cf. [6, p. 62]).  
As is well-known, there is a decomposition of $L^2(G/\Gamma)$ into the orthogonal direct sum of 
\[
L_{\dis}^2 (G/\Gamma) : \qquad \text{the discrete spectrum}
\]
and 
\[
L_{\con}^2 (G/\Gamma) : \qquad \text{the continuous spectrum}.
\]
\vspace{.3cm}

Consider the following statement: 
\vspace{.3cm}
\\
\textbf{Main Conjecture (MC)}.  
The operator $L_{G/\Gamma}^\dis (\alpha)$ is trace class for every \mK-finite $\alpha$ in $C_c^\infty(G)$.
\vspace{.3cm}

This conjecture is a theorem when $\rank (\Gamma ) = 0$ (cf. [6, p. 355])  or when $\rank (\Gamma )= 1$ (cf. Donnelly [3, p. 349]) and is actually true in general (cf. M\"uller [4, p. 473]).

Owing to the theory of the parametrix (cf. [6, p. 21]), it then automatically holds for all \mK-finite $\alpha$ in $\sC^1(G)$.

These points made, the fundamental problem of the theory is to compute
\[
\tr(L_{G/\Gamma}^\dis (\alpha))
\]
in explicit terms.  Thanks to the considerations to be found in [7-(c), \S8], the problem can be divided into two parts:

(1) \quad Determine the contribution to the trace arising from the conjugacy classes.

(2) \quad Determine the contribution to the trace arising from the continuous spectrum.  

[Note: \quad Naturally, when rank $(\Gamma) = 0$, (2) is irrelevant, so only (1) is of interest, an elementary matter.]

Our approach dictates that the second issue be addressed first.  
The essence of the method of attack can be found already in [7-(a)], the key being the cancellation principle.  
There, of course, $\rank (\Gamma) = 1$ and all the Arthur polynomials are linear, so everything, by comparison is fairly simple.  
The situation when $\rank (\Gamma) > 1$ is far more complicated.  
Nevertheless, it is still possible to arrive at an explicit determination, the basis for the cancellation being a certain remarkable ``addition'' property enjoyed by the Arthur polynomials, combined with a multidimensional Dini calculus.  
The way it works is this.  
Each proper \mG-conjugacy class $\sC$ of $\Gamma$-cuspidal split parabolic subgroups of \mG makes a contribution 
\[
\bCon(\alpha : \Gamma : \sC)
\] 
to the trace, the total contribution 
to the trace furnished by the continuous spectrum being the sum
\[
\bCon - \bSp(\alpha : \Gamma) \ = \ \sum_\sC \bCon(\alpha : \Gamma : \sC).
\]
Accordingly, fix a $\sC$ containing 
$P = M \bullet A \bullet N$, say $-$then
\[
\bCon(\alpha : \Gamma : \sC) \ = \ 
\sum\limits_\sO \sum\limits_{w \in W(A)} \ 
\bCon(\alpha : \Gamma : \sC : \sO : w),
\] 
the actual form of the contribution
\[
\bCon(\alpha : \Gamma : \sC : \sO : w)
\] 
depending on $w$ through
\[
\rank (1 - w),
\]
the orbit type $\sO$ having a passive part in the overall procedure.

To provide some motivation for [7-(d)], we shall explicate here the position when $\rank(\Gamma) = 2$.  
Before doing this, though, it will be a good idea to recall how things go when $\rank(\Gamma) = 1$.  
For use below, denote by $*(\sC)$ the number of chambers in \mA (cf. [6, p. 104]).

Fixing $\sO$, let us now suppose that $\rank(\Gamma) = 1$ $-$then $\#(W(A)) = 2$.  
Thus, there are two terms appearing in the contribution from the continuous spectrum.
\vspace{.3cm}

\fbox{$w = 1$} \quad In this case, 
\[
\bCon(\alpha : \Gamma : \sC : \sO : 1)
\] 
is equal to 
\[
-\frac{1}{2 \pi} \bullet \frac{1}{*(\sC)} \bullet  \sum\limits_{w \in W(A)} 
\]
\[
\times \int_{\Re(\Lambda) = 0} \ 
\tr\bigg( \bInd_P^G((\sO,\Lambda))(\alpha)
\]
\[
\bullet \bc ( \restr{P}{A} : \restr{P}{A} : w : \Lambda)^*
\frac{d}{d \Lambda} 
\bc ( \restr{P}{A} : \restr{P}{A} : w : \Lambda)\bigg) \abs{d \Lambda}.
\]

[Note: \quad Since the \bc-function attached to the trivial element of $W(A)$ is a constant, the contribution is concentrated entirely in the \bc-function of the nontrivial element of $W(A)$.  
Still, this mode of expression possesses an inherent symmetry that can be generalized.] 
\vspace{.3cm}

\fbox{$w \neq 1$} \quad In this case, 
\[
\bCon(\alpha : \Gamma : \sC : \sO : w)
\] 
is equal to 
\[
-\frac{1}{2 \pi} \bullet 2 \pi \bullet \frac{1}{*(\sC)} \bullet  \frac{1}{\abs{\det(1 - w)}}
\]
\[
\times \ 
\tr\bigg( \bInd_P^G((\sO,0))(\alpha) \bullet \bc ( \restr{P}{A} : \restr{P}{A} : w : 0) \bigg).
\]

[Note: \quad Since
\[
w \neq 1 \implies \abs{\det(1 - w)} \ = \ 2,
\]
the prefacing constant is 1/4.  
The ``$1/2 \pi$'' is inherent in the Fourier inversion formula; the ``$2 \pi$ is inherent in the Dini calculus.  
Because $1 - w$ is nonsingular, they cancel.]
\vspace{.3cm}

Keeping the orbit type fixed, assume now that $\rank(\Gamma) = 2$.   
There are then two \mG-conjugacy classes $\sC^\prime$ and $\sC^{\prime\prime}$ of maximal $\Gamma$-cuspidal split parabolic subgroups of \mG and one \mG-conjugacy class $\sC$ of minimal $\Gamma$-cuspidal split parabolic subgroups of \mG.  It will be best to discuss each level separately.
\vspace{.3cm}

\fbox{$\sC^\prime, \ \sC^{\prime\prime}$} \quad Two cases can occur.
\vspace{.2cm}


(I) \quad Suppose that $P^\prime \in \sC^\prime$, $P^{\prime\prime} \in \sC^{\prime\prime}$ are associate 
(e.g. $\bA_2$) $-$then $W(A^\prime) = \{1\}$, $W(A^{\prime\prime}) = \{1\}$ and 
\[
\begin{cases}
\ W(A^{\prime\prime},A^\prime)) \ = \ \{w^\prime\}\\
\ W(A^{\prime},A^{\prime\prime})) \ = \ \{w^{\prime\prime}\}
\end{cases}
.
\]
In this case, 
\[
\bCon(\alpha : \Gamma : \sC^\prime : \sO^\prime : 1)
\] 
is equal to 
\[
-\frac{1}{2 \pi} \bullet \frac{1}{*(\sC^\prime)} 
\]
\[
\times \int_{\Re(\Lambda) = 0} \ 
\tr\bigg( \bInd_{P^\prime}^G((\sO^\prime,\Lambda^\prime))(\alpha)
\]
\[
\bullet 
\bc ( \restr{P^{\prime\prime}}{A^{\prime\prime}} : \restr{P^\prime}{A^\prime} : w^\prime : \Lambda^\prime)^*
\frac{d}{d \Lambda^\prime} 
\bc ( \restr{P^{\prime\prime}}{A^{\prime\prime}} : \restr{P^\prime}{A^\prime} : w^\prime : \Lambda^\prime)\bigg) \abs{d \Lambda^\prime}
\]
and 
\[
\bCon(\alpha : \Gamma : \sC^{\prime\prime} : \sO^{\prime\prime} : 1)
\] 
is equal to 
\[
-\frac{1}{2 \pi} \bullet \frac{1}{*(\sC^{\prime\prime})} 
\]
\[
\times \ 
\int_{\Re(\Lambda^{\prime\prime}) = 0} \ 
\tr\bigg( \bInd_{P^{\prime\prime}}^G((\sO^{\prime\prime},\Lambda^{\prime\prime}))(\alpha)
\]
\[
\bullet \bc ( \restr{P^\prime}{A^\prime} : \restr{P^{\prime\prime}}{A^{\prime\prime}} : w^{\prime\prime} : \Lambda^{\prime\prime})^*
\frac{d}{d \Lambda^{\prime\prime}} 
\bc ( \restr{P^\prime}{A^\prime} : \restr{P^{\prime\prime}}{A^{\prime\prime}} : 
w^{\prime\prime} : \Lambda^{\prime\prime})\bigg) \abs{d \Lambda^{\prime\prime}}.
\]
\vspace{.3cm}

(II) \quad Suppose that $P^\prime \in \sC^\prime$, $P^{\prime\prime} \in \sC^{\prime\prime}$ are not associate 
(e.g. $\bA_1 \times \bA_1$, $\bB_2$, $\bG_2$) $-$then 
\[
\begin{cases}
\ W(A^\prime) \ = \ \{1, w^\prime\}\\
\ W(A^{\prime\prime}) \ = \ \{1, w^{\prime\prime}\}
\end{cases}
.
\]
In this case, 
\[
\bCon(\alpha : \Gamma : \sC^\prime : \sO^\prime : 1)
\] 
is equal to 
\[
-\frac{1}{2 \pi} \bullet \frac{1}{*(\sC^\prime)} \bullet  \sum\limits_{w^\prime \in W(A^\prime)} 
\]
\[
\times \int_{\Re(\Lambda^\prime) = 0} \ 
\tr\bigg( \bInd_{P^\prime}^G((\sO^\prime,\Lambda^\prime))(\alpha)
\]
\[
\bullet \bc ( \restr{P^\prime}{A^\prime} : \restr{P^\prime}{A^\prime} : w^\prime : \Lambda^\prime)^*
\frac{d}{d \Lambda^\prime} 
\bc ( \restr{P^{\prime}}{A^{\prime}} : \restr{P^\prime}{A^\prime} : w^\prime : \Lambda^\prime)\bigg) 
\abs{d \Lambda^\prime}
\]
and 
\[
\bCon(\alpha : \Gamma : \sC^{\prime} : \sO^{\prime} : w^\prime)
\] 
is equal to 
\[
-\frac{1}{2 \pi} 
\bullet 
2 \pi \bullet \frac{1}{*(\sC^{\prime})} 
\bullet
\frac{1}{\abs{\det(1 - w^\prime)}}
\]
\[
\times \ 
\tr\bigg( \bInd_{P^\prime}^G((\sO^\prime,0))(\alpha) 
\bullet
\bc (\restr{P^\prime}{A^\prime} : \restr{P^\prime}{A^\prime} : w^\prime : 0)\bigg),
\]
while
\[
\bCon(\alpha : \Gamma : \sC^{\prime\prime} : \sO^{\prime\prime} : 1)
\] 
is equal to 
\[
-\frac{1}{2 \pi} \bullet \frac{1}{*(\sC^{\prime\prime})} \bullet  \sum\limits_{w^{\prime\prime} \in W(A^{\prime\prime})} 
\]
\[
\times \int_{\Re(\Lambda^{\prime\prime}) = 0} \ 
\tr\bigg( \bInd_{P^{\prime\prime}}^G((\sO^{\prime\prime},\Lambda^{\prime\prime}))(\alpha)
\]
\[
\bullet \bc ( \restr{P^{\prime\prime}}{A^{\prime\prime}} 
: \restr{P^{\prime\prime}}{A^{\prime\prime}} : w^{\prime\prime} : \Lambda^{\prime\prime})^*
\frac{d}{d \Lambda^{\prime\prime}} 
\bc ( \restr{P^{\prime\prime}}{A^{\prime\prime}} 
: \restr{P^{\prime\prime}}{A^{\prime\prime}} : w^{\prime\prime} : \Lambda^{\prime\prime})\bigg) 
\abs{d \Lambda^{\prime\prime}}  
\]
and 
\[
\bCon(\alpha : \Gamma : \sC^{\prime\prime} : \sO^{\prime\prime} : w^{\prime\prime})
\] 
is equal to 
\[
-\frac{1}{2 \pi} \bullet 2 \pi \bullet \frac{1}{*(\sC^{\prime\prime})} \frac{1}{\abs{\det(1 - w^{\prime\prime})}}
\]
\[
\times \ 
\tr\bigg( \bInd_{P^{\prime\prime}}^G((\sO^{\prime\prime},0))(\alpha) 
\bullet \bc ( \restr{P^{\prime\prime}}{A^{\prime\prime}} : \restr{P^{\prime\prime}}{A^{\prime\prime}} :
w^{\prime\prime} :0) \bigg).
\]
\vspace{.05cm}


[Note: \quad At the maximal level, therefore, the contribution to the trace is entirely analogous to what obtains when 
$\rank(\Gamma) = 1$, including the interpretation of the constants.]
\vspace{.3cm}

\fbox{$\sC$} \quad 
Given $w \in W(A)$, there are three possibilities:
\[
\begin{cases}
\ \rank(1 - w) = 0\\
\ \rank(1 - w) = 1\\
\ \rank(1 - w) = 2
\end{cases}
.
\]
The two extreme cases are the easiest to treat and will be dealt with first.

Let $\lambda_1$ and $\lambda_2$ be the simple roots: let $\lambda^1$ and $\lambda^2$ be their duals.  
Generically, write
\[
\widehat{\lambda} \ = \ \frac{\lambda}{\norm{\lambda}}.
\]
\vspace{.3cm}

\fbox{$\rank (1 - w) = 0$} \quad 
This requirement implies that $w = 1$.   Introduce
\[
P_P^G(H) \ = \ \frac{1}{2}
\big\{(H,\widehat{\lambda}_1)(H,\widehat{\lambda}^2) + (H,\widehat{\lambda}^1)(H,\widehat{\lambda}_2)\big\}.
\]
Then $P_P^G$ is an Arthur polynomial.  
As such, it is homogeneous of degree 2.  
Denote by $D_P^G$ the associated differential operator.  
In this case, 
\[
\bCon(\alpha : \Gamma : \sC : \sO : 1)
\] 
is equal to 
\[
-\frac{1}{(2 \pi)^2} \bullet \frac{1}{*(\sC)} \bullet  \sum\limits_{w \in W(A)} 
\]
\[
\times \int_{\Re(\Lambda) = 0} \ 
\tr\bigg( \bInd_{P}^G((\sO,\Lambda))(\alpha)
\]
\[
\bullet \bc ( \restr{P}{A} : \restr{P}{A} : w : \Lambda)^*
D_P^G 
\bc ( \restr{P}{A} : \restr{P}{A} : w : \Lambda)\bigg) 
\abs{d \Lambda}.
\]
\vspace{.2cm}

[Note: \quad 
The similarity with the ``$w = 1$'' contribution when $\rank(\Gamma) = 1$ is quite striking.  
In particular, the constants have the ``right'' interpretation and the derivative is ``logarithmic'' in character.  
Needless to say, in the sum over $w \in W(A)$, the term corresponding to $w = 1$ is a priori, zero.]
\vspace{.3cm}

\fbox{$\rank (1 - w) = 2$} \quad
The requirement implies that $1 - w$ is nonsingular.  
In this case, 
\[
\bCon(\alpha : \Gamma : \sC : \sO : w)
\]
is equal to 
\[
-\frac{1}{(2 \pi)^2} 
\bullet
(2 \pi)^2 
\bullet
\frac{1}{*(\sC)} 
\bullet
\frac{1}{\abs{\det(1 - w)}} 
\]
\[
\times \ 
\tr \bigg( \bInd_P^G ((\sO,0))(\alpha) \bullet \bc(\restr{P}{A} : \restr{P}{A} : w : 0) \bigg).
\]
[Note: \quad 
Again, the resemblance to the ``$w \neq 1$'' contribution when $\rank(\Gamma) = 1$ is immediately apparent.  
Once more, the ``$1/(2 \pi)^2$'' is inherent in the Fourier inversion formula; the $(2 \pi)^2$ in inherenet in the Dini calculus.  
Because $1 - w$ is nonsingular, they cancel.]
\vspace{.3cm}

\fbox{$\rank (1 - w) = 1$} \quad
This requirement implies that $w$ is a reflection, say $w = w_\lambda$, where, without loss of generality, $\lambda$ is a positive short root.  
The extra ``$1/2$'' that arises in what follows has its origins in a change of variables, which can be traced back to the fact that 
\[
(1 - w_\lambda)(\lambda) \ = \ 2\lambda.
\]
Put
\[
*(\sC(\lambda)) \ = \ 2
\]
and let $D_\lambda$ be the differential operator corresponding to $-\widehat{\lambda}$.  
We distinguish two cases.  
\vspace{.3cm}


(I) \quad
$\lambda \not\perp \lambda_1$ and $\lambda \not\perp \lambda_2$.  
Write
\[
\begin{cases}
\ \theta_1 \quad \text{for the angle between $\lambda$ and $\lambda_1$}\\
\ \theta_2 \quad \text{for the angle between $\lambda$ and $\lambda_2$}
\end{cases}
.
\]
In this case, 
\[
\bCon(\alpha : \Gamma : \sC : \sO : w_\lambda)
\]
is equal to
\[
-\frac{1}{(2 \pi)^2} 
\bullet
2 \pi
\bullet
\frac{1}{* (\sC(\lambda))}
\bullet
\frac{1}{\abs{\restr{\det ((1 - w_\lambda)}{\Ker(1 - w_\lambda)^\perp)}}}
\bullet
\frac{1}{2}
\]
\[
\bullet 
\frac{sin(\theta_1 + \theta_2)}{cos(\theta_1) cos(\theta_2)}
\]
\[
\times \ 
D_\lambda |_{\Lambda^\prime = 0} \ 
\int_{\Ker(1 - w_\lambda)} \tr\bigg( \bInd_P^G((\sO, \Lambda + \Lambda^\prime))(\alpha)
\]
\[
\bullet \bc(\restr{P}{A} : \restr{P}{A} : w_\lambda : \Lambda + \Lambda^\prime)\bigg) \abs{d \Lambda}.
\]
\vspace{.1cm}

[Note: \quad 
The constant
\[
\frac{sin(\theta_1 + \theta_2)}{cos(\theta_1) cos(\theta_2)}
\]
is strictly positive or strictly negative.]
\vspace{.3cm}

(II) \quad 
$\lambda \perp \lambda_1$ or $\lambda \perp \lambda_2$.  
Let $i = 1$ or 2 and suppose that $\lambda \perp \lambda_i$.  
In this case
\[
\bCon(\alpha : \Gamma : \sC : \sO : w_\lambda)
\]
is equal to the sum of a pair of terms, namely:
\vspace{.3cm}

(II$_1$) \quad 
Call $w_i$ the simple reflection in $\lambda_i$ $-$then the first term is 
\[
-\frac{1}{(2 \pi)^2} 
\bullet
2 \pi
\bullet
\frac{1}{* (\sC(\lambda))}
\bullet
\frac{1}{\abs{\restr{\det ((1 - w_\lambda)}{\Ker(1 - w_\lambda)^\perp)}}}
\bullet
\frac{1}{2}
\]
\[
\times \ \int_{\Ker (1 - w_\lambda)} \tr \bigg( \bInd_P^G((\sO,\Lambda))(\alpha)
\]
\[
\bullet \bc(\restr{P}{A} : \restr{P}{A} : w_i w_\lambda : \Lambda)^* 
\frac{d}{d \Lambda} \bc(\restr{P}{A} : \restr{P}{A} : w_i : \Lambda)\bigg)
\abs{d \Lambda}.
\]
\vspace{.05cm}

[Note: \quad Here, the \bc-function enters as a ``hybrid'' logarithmic derivative.]

\vspace{.3cm}

(II$_2$) \quad 
Call $\theta_{1 2}$ the angle between $\lambda_1$ and $\lambda_2$ $-$then the second term is 
\[
-\frac{1}{(2 \pi)^2} 
\bullet
2 \pi
\bullet
\frac{1}{* (\sC(\lambda))}
\bullet
\frac{1}{\abs{\restr{\det ((1 - w_\lambda)}{\Ker(1 - w_\lambda)^\perp)}}}
\bullet
\frac{1}{2}
\bullet
\cot(\pi - \theta_{12})
\]
\[
\times \ 
D_\lambda |_{\Lambda^\prime = 0}
\int_{\Ker (1 - w_\lambda)} \tr \bigg( \bInd_P^G((\sO,\Lambda + \Lambda^\prime))(\alpha)
\]
\[
\bullet 
\bc(\restr{P}{A} : \restr{P}{A} : w_\lambda : \Lambda + \Lambda^\prime)\bigg)
\abs{d \Lambda}.
\]

[Note: \quad Since $(\lambda_1, \lambda_2)$ is $\leq 0$, the cotangent of $\pi - \theta_{1 2}$ is $\geq 0$ and can $= 0$ 
(e.g. in $\bA_1 \times \bA_2)$.]
\vspace{.2cm}

We remark that the ``$2 \pi$'' supra is the Dirac constant, hence does not cancel the ``$1 / (2 \pi)^2$'', the Fourier constant.  
Also, $\forall \ \lambda$, 
\[
\abs{\restr{\det(( 1 - w_\lambda)}{\Ker(1 - w_\lambda)^\perp)}} \ = \ 2.
\]

To have a specific illustration of all this, take
\[
\begin{cases}
\ G = \bSL(3,\R)\\
\ \Gamma = \bSL(3,\Z)
\end{cases}
.
\]
Then $\#(W(A)) = 6$.  
Apart from $w = 1$, there are two rotations, $w^\prime$ and $w^{\prime\prime}$, and three reflections, $w_1$, $w_2$, and $w_3$.  
Regarding the latter, only case I applies and we accordingly pick up a sum
\[
\sum\limits_{i = 1}^3 \ \bCon(\alpha : \Gamma : \sC : \sO : w_i)
\]
of three ``orthogonal derivatives''.  

The appearance of 
\[
D_\lambda |_{\Lambda^\prime = 0} \ 
\int_{\Ker(1 - w_\lambda)} \tr\bigg( \bInd_P^G((\sO, \Lambda + \Lambda^\prime))(\alpha)
\]
\[
\bullet \bc(\restr{P}{A} : \restr{P}{A} : w_\lambda : \Lambda + \Lambda^\prime)\bigg) \abs{d \Lambda}
\]
is not a total surprise, if only because in higher rank derivatives of Dirac distributions are produced by the Dini calculus in the presence of quadratic denominators (via the two roots).  
Indeed, if 
\[
\delta (\Ker(1 - w_\lambda))
\]
is the Dirac distribution concentrated on $\Ker(1 - w_\lambda)$, then our ``orthogonal derivative'' is, up to a constant, the result of applying 
\[
\delta^\prime (\Ker(1 - w_\lambda))
\]
to
\[
\tr\bigg( \bInd_P^G((\sO,?))(\alpha) \bullet \bc(\restr{P}{A} : \restr{P}{A} : w_\lambda : ?)\bigg).
\]
\vspace{.3cm}



\begingroup
\center {\textbf{REFERENCES}}\\
\endgroup
\vspace{0.5cm}

\noindent Arthur, J.:

[1] \quad 
A trace formula for $\bGL(3)$, Preprint.


\noindent Arthur, J.: 

[2-(a)] \quad
{On a family of distributions obtained from Eisenstein series I}, 
\emph{Amer. J. Math.}
\textbf{104} (1982), 
1243-1288.

[2-(b)] \quad
{On a family of distributions obtained from Eisenstein series II}, 
\emph{Amer. J. Math.}
\textbf{104} (1982), 
1289-1336.


\noindent Donnelly, H.: 

[3] \quad
{Eigenvalue estimates for certain noncompact manifolds}, 
\emph{Michigan. Math. J.}
\textbf{31} (1984), 
349-357.





\noindent M\"uller, W.: 

[4] \quad
{The trace class conjecture in the theory of automorphic forms}, 
\emph{Ann. of Math.}
\textbf{130} (1989), 
473-529.

\noindent Osborne, S.: 

[5] \quad
{The continuous spectrum: Some explicit formulas, Bowdoin Conference on the Trace Formula}, 
\emph{Contem. Math.},
\textbf{53} (1986), 
371-373.

\begingroup
\noindent {Osborne, M. S., and Warner, G.:}
\endgroup

[6] \quad
{The theory of Eisenstein systems}, 
\emph{Academic Press}, 
New York (1981).


\begingroup
\noindent {Osborne, M. S., and Warner, G.}:
\endgroup

[7-(a)] \quad
{The Selberg trace formula}, 
\emph{Crelle's J.}
\textbf{324} (1981).  
1-113.

[7-(b)] \quad
{The Selberg trace formula V}, 
\emph{Trans. Amer. Math. Soc.}
\textbf{286} (1984), 
351-376.

[7-(c)] \quad
{The Selberg trace formula VI}, 
\emph{Amer. J. Math.}
\textbf{107} (1985), 
1369-1437.

[7-(d)] \quad
{The Selberg trace formula VII}, 
\emph{Pacific J. Math.}
\textbf{140} (1989), 
263-352.

[7-(e)] \quad
{The Selberg trace formula VIII}, 
\emph{Trans. Amer. Math. Soc.}
\textbf{324} (1991), 
623-653.

[7-(f)] \quad
{The Selberg trace formula IX}, 
Preprint.

\begingroup
\noindent {Osborne, M. S., and Warner, G.:}
\endgroup

[8] \quad
{On the Dini calculus}, 
\emph{J. Math. Anal. and Applic.}
\textbf{118} (1986), 
255-262 (see also \textbf{139} (1989), 282-300).


\begingroup
\noindent {Venkov, A.:}
\endgroup

[9] \quad
{The Selberg trace formula for $\bSL(3,\Z)$}, 
\emph{J. Soviet Math.}
\textbf{12} (1979), 
384-424.

\noindent Warner, G.: 

[10-(a)] \quad
{Selberg's trace formula for nonuniform lattices: The $\R$-rank 1 case}, 
\emph{Adv. in Math. Studies}
\textbf{6} (1979), 
1-142.

[10-(b)] \quad
{Traceability on the discrete spectrum, Bowdoin Conference on the Trace Formula}, 
\emph{Contem. Math.},
\textbf{53} (1986), 
529-534.


\end{document}